\documentclass[12pt]{article}%
\usepackage{amsmath}
\usepackage{amsfonts}
\usepackage{amssymb}
\usepackage{graphicx}%
\providecommand{\U}[1]{\protect\rule{.1in}{.1in}}
\newtheorem{theorem}{Theorem}[section]

\newtheorem{claim}[theorem]{Claim}

\newtheorem{corollary}[theorem]{Corollary}

\newtheorem{lemma}[theorem]{Lemma}

\newtheorem{observation}[theorem]{Observation}
\newtheorem{problem}{Problem}
\newtheorem{proposition}[theorem]{Proposition}
\newtheorem{question}{Question}

\begin{document}

\title{\textbf{Lower Bound and Exact Values for the Boundary Independence Broadcast
Number of a Tree}}
\author{C.M. Mynhardt\thanks{Supported by the Natural Sciences and Engineering
Research Council of Canada, PIN 253271.}\\Department of Mathematics and Statistics\\University of Victoria, Victoria, BC, \textsc{Canada}\\{\small kieka@uvic.ca}
\and L. Neilson\\Department of Adult Basic Education\\Vancouver Island University, Nanaimo, BC, \textsc{Canada}\\{\small linda.neilson@viu.ca}}
\maketitle

\begin{abstract}
A broadcast on a nontrivial connected graph $G=(V,E)$ is a function
$f:V\rightarrow\{0,1,\dots,\operatorname{diam}(G)\}$ such that $f(v)\leq e(v)$
(the eccentricity of $v$) for all $v\in V$. The weight of $f$ is $\sigma(f)=%
{\textstyle\sum_{v\in V}}
f(v)$. A vertex $u$ hears $f$ from $v$ if $f(v)>0$ and $d(u,v)\leq f(v)$.

A broadcast $f$ is boundary independent\emph{ }if, for any vertex $w$ that
hears $f$ from vertices $v_{1},...,v_{k},\ k\geq2$, we have that
$d(w,v_{i})=f(v_{i})$ for each $i$. The maximum weight of a boundary
independent broadcast on $G$ is denoted by\emph{ }$\alpha_{\operatorname{bn}%
}(G)$. We prove a sharp lower bound on $\alpha_{\operatorname{bn}}(T)$ for a
tree $T$. Combined with a previously determined upper bound, this gives exact
values of $\alpha_{\operatorname{bn}}(T)$ for some classes of trees $T$. We
also determine $\alpha_{\operatorname{bn}}(T)$ for trees with exactly two
branch vertices and use this result to demonstrate the existence of trees for
which $\alpha_{\operatorname{bn}}$ lies strictly between the lower and upper bounds.

\end{abstract}

\noindent\textbf{Keywords:\hspace{0.1in}}broadcast domination; broadcast
independence; hearing independence; boundary independence

\noindent\textbf{AMS Subject Classification Number 2010:\hspace{0.1in}}05C69

\section{Introduction}

An independent set $X$ of vertices in a graph $G$, when considered from the
perspective of its vertices, has the property that no vertex in $X$ belongs to
the neighbourhood of any other vertex in $X$ (i.e., vertices in $X$ are
nonadjacent). Seen from the perspective of the edges of $G$, it has the
property that no edge is covered by more than one vertex in $X$, that is, the
neighbourhoods of vertices in $X$ overlap only on their boundaries. In
generalizing independent sets to independent broadcasts, Erwin \cite{Ethesis}
used the former property, and defined a broadcast to be independent if no
broadcasting vertex belongs to the neighbourhood of, or hears, another
broadcasting vertex. We refer to this type of independent broadcast as being
hearing independent. Focussing instead on the latter property, Mynhardt and
Neilson \cite{MN} and Neilson \cite{LindaD} defined boundary independent
broadcasts as broadcasts in which any edge is covered by at most one
broadcasting vertex, that is, a vertex belongs to the broadcasting
neighbourhoods of two or more broadcasting vertices only if it belongs to the
boundaries of all such vertices.

Our goal here is to prove Theorem \ref{Thm_LB}, which gives a tight lower
bound on the boundary independence number of a tree. Together with an upper
bound proved by Mynhardt and Neilson \cite{MN3}, this gives exact values of
the boundary independence numbers of some classes of trees. We also determine
an exact formula in Theorem \ref{b(T)2} for boundary independence numbers of
trees with exactly two branch vertices, and use this result to show that there
exist trees whose boundary independence numbers lie strictly between the
above-mentioned lower and upper bounds. 

After giving the necessary definitions in Sections \ref{Sec_broadcast_Defs}
and \ref{SecTreeDefs}, we state the lower bound and the formula in Section
\ref{SecMainThm}. The proofs are given in Sections \ref{Sec_lower_pf} and
\ref{Thm_b(T)2_pf}, respectively. For comparison with the lower bound we also
give the upper bound from \cite{MN3} in Section~\ref{SecMainThm}, and discuss
some corollaries of the bounds. Section \ref{Sec_known} contains previous
results and lemmas required for the proofs of the main results. We conclude by
listing open problems in Section \ref{Sec_Open}. For undefined concepts we
refer the reader to~\cite{CLZ}.

\subsection{Broadcast definitions}

\label{Sec_broadcast_Defs}As this paper is a sequel to the paper \cite{MN3} by
the same authors, the definitions and background information provided in this
and the next subsections are essentially the same as those in \cite{MN3}. A
\emph{broadcast} on a nontrivial connected graph $G=(V,E)$, as introduced by
Erwin \cite{Ethesis, Epaper}, is a function $f:V\rightarrow\{0,1,\dots
,\operatorname{diam}(G)\}$ such that $f(v)\leq e(v)$ (the eccentricity of $v$)
for all $v\in V$. If $G$ is disconnected, a broadcast on $G$ is the union of
broadcasts on its components. The \emph{weight} of $f$ is $\sigma
(f)=\sum_{v\in V}f(v)$. Define $V_{f}^{+}=\{v\in V:f(v)>0\}$ and partition
$V_{f}^{+}$ into the two sets $V_{f}^{1}=\{v\in V:f(v)=1\}$ and $V_{f}%
^{++}=V_{f}^{+}-V_{f}^{1}$. A vertex in $V_{f}^{+}$ is called a
\emph{broadcasting vertex}. A vertex $u$ \emph{hears} $f$ from $v\in V_{f}%
^{+}$, and $v$ $f$-\emph{dominates} $u$, if the distance $d(u,v)\leq f(v)$. If
$d(u,v)<f(v)$, we also say that say that $v$ \emph{overdominates }$u$. Denote
the set of all vertices that do not hear $f$ by $U_{f}$. A broadcast $f$ is
\emph{dominating} if $U_{f}=\varnothing$. If $f$ is a broadcast such that
every vertex $x$ that hears more than one broadcasting vertex also satisfies
$d(x,u)\geq f(u)$ for all $u\in V_{f}^{+}$, we say that the \emph{broadcast
only overlaps in boundaries}. If $uv\in E(G)$ and $u,v\in N_{f}(x)$ for some
$x\in V_{f}^{+}$ such that at least one of $u$ and $v$ does not belong to
$B_{f}(x)$, we say that the edge $uv$ is \emph{covered} in $f$, or
$f$-\emph{covered}, by $x$. If $uv$ is not covered by any $x\in V_{f}^{+}$, we
say that $uv$ is \emph{uncovered by~}$f$ or $f$-\emph{uncovered}. We denote
the set of $f$-uncovered edges by $U_{f}^{E}$.

If $f$ and $g$ are broadcasts on $G$ such that $g(v)\leq f(v)$ for each $v\in
V$, we write $g\leq f$. If in addition $g(v)<f(v)$ for at least one $v\in V$,
we write $g<f$. A dominating broadcast $f$ on $G$ is a \emph{minimal
dominating broadcast} if no broadcast $g<f$ is dominating. The \emph{upper
broadcast number }of $G$ is $\Gamma_{b}(G)=\max\left\{  \sigma(f):f\text{ is a
minimal dominating broadcast of }G\right\}  $. First defined by Erwin
\cite{Ethesis}, the upper broadcast number was also studied by, for example,
Ahmadi, Fricke, Schroeder, Hedetniemi and Laskar \cite{Ahmadi}, Bouchemakh and
Fergani \cite{BF}, Bouchouika, Bouchemakh and Sopena \cite{BBS}, Dunbar,
Erwin, Haynes, Hedetniemi and Hedetniemi \cite{DEHHH}, and Mynhardt and Roux
\cite{MR}.

We denote the independence number of $G$ by $\alpha(G)$. If $f$ is
characteristic function of an independent set of $G$, then no vertex in
$V_{f}^{+}$ hears $f$ from any other vertex. To generalize the concept of
independent sets, Erwin \cite{Ethesis} defined a broadcast $f$ to be
\emph{independent}, or, for our purposes, \emph{hearing independent}, if no
vertex $u\in V_{f}^{+}$ hears $f$ from any other vertex $v\in V_{f}^{+}$; that
is, broadcasting vertices only hear themselves. We denote the maximum weight
of a hearing independent broadcast on $G$ by $\alpha_{h}(G)$. This version of
broadcast independence was also considered by, among others, Ahmane,
Bouchemakh and Sopena \cite{ABS, ABS2}, Bessy and Rautenbach \cite{BR, BR2},
Bouchemakh and Zemir \cite{Bouch}, Bouchouika et al.~\cite{BBS} and Dunbar et
al.~\cite{DEHHH}. For a survey of broadcasts in graphs, see the chapter by
Henning, MacGillivray and Yang \cite{HMY}.

For a broadcast $f$ on a graph $G$ and $v\in V_{f}^{+}$, we define the%
\[
\left.
\begin{tabular}
[c]{l}%
$f$-\emph{neighbourhood}\\
$f$-\emph{boundary}\\
$f$-\emph{private neighbourhood}\\%
\begin{tabular}
[c]{l}%
\
\end{tabular}
\\
$f$-\emph{private boundary}\\%
\begin{tabular}
[c]{l}%
\
\end{tabular}
\end{tabular}
\ \ \ \right\}  \ \text{of }v\text{ by }\left\{  \text{%
\begin{tabular}
[c]{rll}%
$N_{f}(v)$ & $=$ & $\{u\in V:d(u,v)\leq f(v)\}$\\
$B_{f}(v)$ & $=$ & $\{u\in V:d(u,v)=f(v)\}$\\
$\operatorname{PN}_{f}(v)$ & $=$ & $\{u\in N_{f}(v):u\notin N_{f}(w)$ for
all\\
&  & $\ w\in V_{f}^{+}-\{v\}\}$\\
$\operatorname{PB}_{f}(v)$ & $=$ & $\{u\in N_{f}(v):u$ is not dominated by\\
&  & $\ (f-\{(v,f(v)\})\cup\{(v,f(v)-1)\}$.
\end{tabular}
}\right.
\]

If $v\in V_{f}^{1}$ and $v$ does not hear $f$ from any vertex $u\in V_{f}%
^{+}-\{v\}$, then $v\in\operatorname{PB}_{f}(v)$, and if $v\in V_{f}^{++}$,
then $\operatorname{PB}_{f}(v)=B_{f}(v)\cap\operatorname{PN}_{f}(v)$. Also
note that $f$ is a broadcast that overlaps only in boundaries if and only if
$N_{f}(u)\cap N_{f}(v)\subseteq B_{f}(u)\cap B_{f}(v)$ for all distinct
$u,v\in V_{f}^{+}$.

The characteristic function of an independent set also has the feature that it
only overlaps in boundaries. To generalize this property, we define a
broadcast to be \emph{boundary independent}, abbreviated to
\emph{bn-independent}, if it overlaps only in boundaries. The maximum weight
of a bn-independent broadcast on $G$ is the \emph{boundary independence
number} $\alpha_{\operatorname{bn}}(G)$; such a broadcast is called an
$\alpha_{\operatorname{bn}}(G)$\emph{-broadcast}, often abbreviated to
$\alpha_{\operatorname{bn}}$\emph{-broadcast} if the graph $G$ is clear. The
respective definitions imply that $\alpha(G)\leq\alpha_{\operatorname{bn}%
}(G)\leq\alpha_{h}(G)$ for all graphs $G$. Boundary independent broadcasts
were introduced by Neilson \cite{LindaD} and Mynhardt and Neilson \cite{MN},
and also studied in \cite{MM, MN2, MN3}. For example, it was shown in
\cite{MN} that $\alpha_{h}(G)/\alpha_{\operatorname{bn}}(G)<2$ for all graphs
$G$, and the bound is asymptotically best possible. In \cite{MN2} it was shown
that $\alpha_{\operatorname{bn}}$ and $\Gamma_{b}$ are not comparable, and
that $\alpha_{\operatorname{bn}}(G)/\Gamma_{b}(G)<2$ for all graphs $G$, while
$\Gamma_{b}(G)/\alpha_{\operatorname{bn}}(G)$ is unbounded. A sharp upper
bound for $\alpha_{\operatorname{bn}}(T)$, where $T$ is a tree, was proved in
\cite{MN3} (see Theorem \ref{ThmMain}).

\subsection{Definitions for trees}

\label{SecTreeDefs}The statements of the upper bound in \cite{MN3} for the
boundary independence number of trees, and the lower bound in this paper,
require some definitions of concepts pertaining to trees. A vertex of a tree
$T$ of degree $3$ or more is called a \emph{branch vertex}. We denote the set
of leaves of $T$ by $L(T)$, the set of branch vertices by $B(T)$ and the set
of vertices of degree $2$ by $W(T)$. The unique neighbour of a leaf is called
a \emph{stem}.

For integers $k\geq3$ and $n_{i}\geq1,\ i\in\{1,...,k\}$, the
\emph{(generalized) spider} $\operatorname{Sp}(n_{1},...,n_{k})$ is the tree
which has exactly one branch vertex $b$, called the \emph{head}, with
$\deg(b)=k$, and for which the $k$ components of $\operatorname{Sp}%
(n_{1},...,n_{k})-b$ are paths of lengths $n_{1}-1,...,n_{k}-1$, respectively.
The \emph{legs }$L_{1},...,L_{k}$ of the spider are the paths from $b$ to the
leaves. Let $t_{i}$ be the leaf of $L_{i},\ i=1,...,k$. If $n_{i}=r$ for each
$i$, we write $\operatorname{Sp}(n_{1},...,n_{k})=\operatorname{Sp}(r^{k})$. A
\emph{caterpillar} of length $k\geq0$ is a tree such that removing all leaves
produces a path of length $k$, called the \emph{spine}.

\begin{figure}[ptb]%
\centering
\includegraphics[
height=3.0217in,
width=5.6265in
]%
{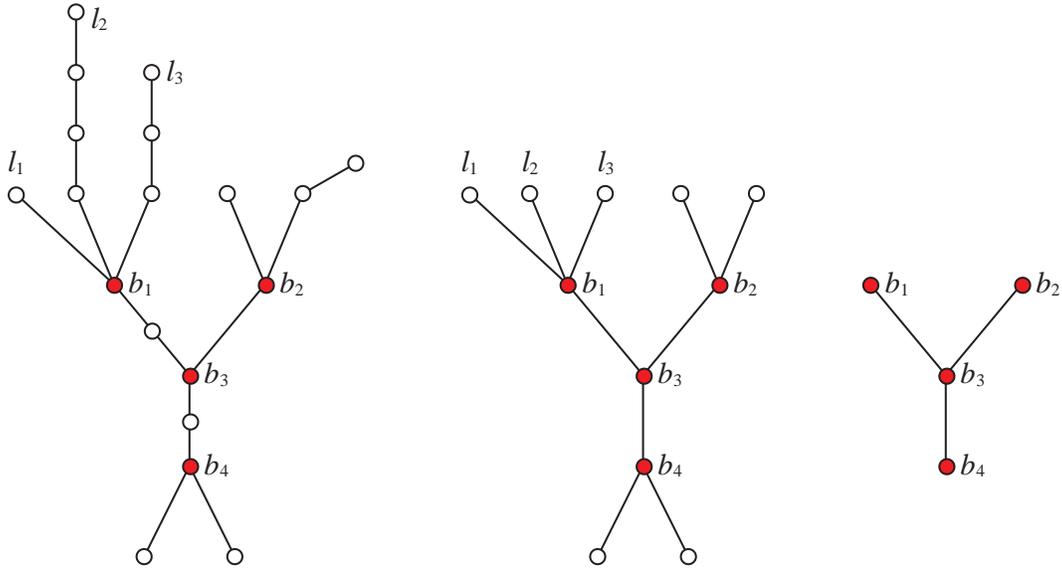}%
\caption{\cite{MN3, LindaD}\ \ A tree $T$ (left), its branch-leaf
representation $\mathcal{BL}(T)$ (middle) and its branch representation
$\mathcal{B}(T)$ (right). The branch set of $T$ is $B(T)=\{b_{1},b_{2}%
,b_{3},b_{4}\}$, $R(T)=\{b_{3}\}$, $B_{\operatorname{end}}(T)=\{b_{1}%
,b_{2},b_{4}\}$, $b(T)=4$ and $\rho(T)=1$. The branch vertex $b_{1}$ has leaf
set $L(b_{1})=\{l_{1},l_{2},l_{3}\}$. }%
\label{Fig6.1}%
\end{figure}

The next few concepts are illustrated in Figure \ref{Fig6.1}. The
\emph{branch-leaf representation }$\mathcal{BL}(T)$ of $T$ is the tree
obtained by suppressing all vertices $v$ with $\deg(v)=2$, and the
\emph{branch representation }$\mathcal{B}(T)$ of a tree $T$ with at least one
branch vertex is obtained by deleting all leaves of $\mathcal{BL}(T)$. Thus,
$V(\mathcal{B}(T))=B(T)$, and two vertices $b_{1},b_{2}\in B(T)$ are adjacent
in $\mathcal{B}(T)$ if and only if the $b_{1}-b_{2}$ path in $T$ contains no
other branch vertices. We denote $|B(T)|$ by $b(T)$.%

An \emph{endpath} in a tree is a path ending in a leaf and having all internal
vertices (if any) of degree $2$. If there exists a $v-l$ endpath, where $v\in
B(T)$ and $l\in L(T)$, then $v$ and $l$ are adjacent \underline{in
$\mathcal{BL}(T)$}; we also say that $l$ belongs to $L(v)$, the \emph{leaf
set} of $v$, and we refer to $l$ as \emph{a leaf of} $v$ (even though $l$ is
not necessarily adjacent to $v$ in $T$). Since $\mathcal{BL}(T)$ is unique, we
can talk about $L(v)$ for any branch vertex $v$, where the reference to
$\mathcal{BL}(T)$ is implied but not specifically mentioned. Let $R(T)$ be the
set of all branch vertices $w$ of $T$ such that $|L(w)|\leq1$ and define
$\rho(T)=|R(T)|$. Equivalently, $\rho(T)$ is the number of branch vertices of
$T$ with at most one leaf, that is, the branch vertices which belong to at
most one endpath. 

\begin{figure}[ptb]%
\centering
\includegraphics[
height=1.8974in,
width=6.3218in
]%
{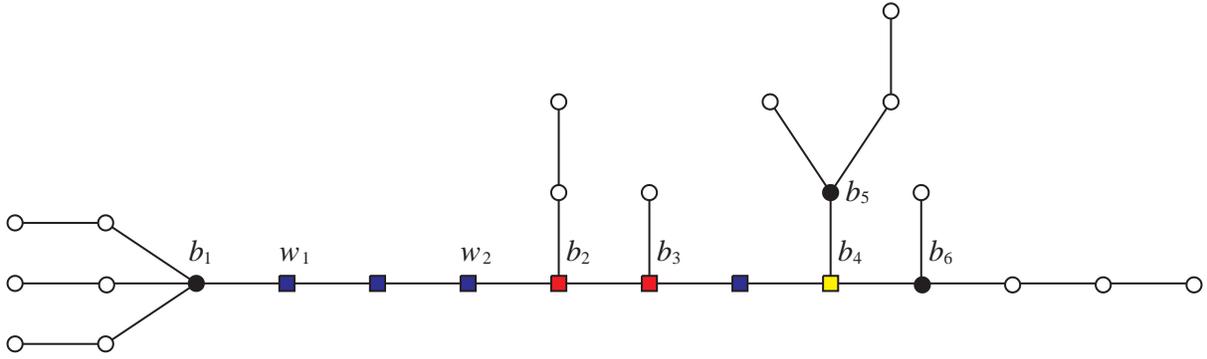}%
\caption{A tree $T$ with the vertices of $W_{\operatorname{int}}(T)$ coloured
blue, $B_{0}(T)$ coloured yellow, $B_{1}(T)$ coloured red, $B_{\geq2}(T)$
coloured black, and $W_{\operatorname{ext}}(T)$ and the leaves uncoloured. The
subgraph $\operatorname{Int}(T)$ is induced by the square vertices;
$\{w_{1},w_{2},b_{3},b_{4}\}$ is a maximal independent set of
$\operatorname{Int}(T)$, and $R(T)=\{b_{2},b_{3},b_{4}\}$.}%
\label{Fig_LB1}%
\end{figure}

We define subsets $B_{i}(T)$ and $B_{\geq i}(T)$ of $B(T)$ by
\[
B_{i}(T)=\{v\in B(T):|L(v)|=i\}\text{ and }B_{\geq i}(T)=\{v\in
B(T):|L(v)|\geq i\}.
\]
Clearly, $B_{0}(T)\cup B_{1}(T)\cup B_{\geq2}(T)$ is a partition of $B(T)$
while $B_{0}(T)\cup B_{1}(T)$ is a partition of $R(T)$. See Figure
\ref{Fig_LB1}. We also partition the set $W(T)$ of vertices of degree $2$ into
two subsets, $W_{\operatorname{ext}}(T)$ for the \emph{external vertices of
degree }$2$, and $W_{\operatorname{int}}(T)$ for the \emph{internal vertices
of degree }$2$, as follows:%
\[
W_{\operatorname{ext}}(T)=\{u\in W(T):u\ \text{lies\ on\ an\ endpath}%
\}\text{\ and\ }W_{\operatorname{int}}(T)=W(T)-W_{\operatorname{ext}}(T).
\]
The subgraph of $T$ induced by $B_{0}(T)\cup B_{1}(T)\cup
W_{\operatorname{int}}(T)$ is called the \emph{interior subgraph} of $T$,
denoted by $\operatorname{Int}(T)$. Note that $\operatorname{Int}(T)$ is not
necessarily connected.%

For a branch vertex $b$, let $T_{b}$ be the subtree of $T$ induced by all the
$b-l$ paths from $b$ to leaves $l\in L(b)$. Then $T_{b}=K_{1}$ if $b\in
B_{0}(T)$, $T_{b}$ is a path of length $d(l,b)$ if $b\in B_{1}(T)$ and
$L(b)=\{l\}$, $T_{b}$ is a path of length $d(l_{1},b)+d(l_{2},b)$ if
$L(b)=\{l_{1},l_{2}\}$, and $T_{b}$ is a generalized spider $S(d(l_{1}%
,b),,...,d(l_{k},b))$ if $L(b)=\{l_{1},...,l_{k}\},\ k\geq3$. Observe that
$\operatorname{Int}(T)$ can also be obtained by deleting the vertices of
$T_{b}$ for each $b\in B_{\geq2}(T)$, and the $b-l$ path, except for $b$, if
$b\in B_{1}(T)$ and $L(b)=\{l\}$. For a tree $T$ with $b(T)\geq1$, the
\label{maxleaf}\label{sumleaf}\label{lossleaf}%
\[
\left.
\begin{tabular}
[c]{r}%
{\emph{max leaf value}}\\
{\emph{sum}}\\
{\emph{loss}}%
\end{tabular}
\ \right\}  \text{ of a branch vertex }v\text{ is\ }\left\{
\begin{tabular}
[c]{l}%
{$\max(v)=\max\{d_{T}(v,x):x\in L(v)\},$}\\
$\operatorname{sum}(v)=\sum_{x\in L(v)}${$d_{T}(v,x),$}\\
{$\operatorname{loss}(v)=\operatorname{sum}(v)-\max(v).$}%
\end{tabular}
\ \ \right.
\]
For the tree $T$ in Figure \ref{Fig6.1}, $\operatorname{sum}(b_{1})=8$,
$\max(b_{1})=4$ and $\operatorname{loss}(b_{1})=4$.

\section{Statement of main results and their corollaries}

\label{SecMainThm}With the necessary definitions in place, we now state the
lower bound for $\alpha_{\operatorname{bn}}(T)$, where $T$ is a tree, as well
as the upper bound in \cite{MN3} for comparison. We defer the proof of the
lower bound to Section \ref{Sec_lower_pf}, after stating some known results
required for the proof.

\begin{theorem}
\label{Thm_LB}If $T$ is a tree such that $b(T)\geq1$, then%
\[
\alpha_{\operatorname{bn}}(T)\geq n-b(T)-|W_{\operatorname{int}}%
(T)|+\alpha(\operatorname{Int}(T)).
\]

\end{theorem}

\begin{theorem}
\label{ThmMain}\emph{\cite[Theorem 1.1]{MN3}}\hspace{0.1in}For any tree $T$ of
order $n$, $\alpha_{\operatorname{bn}}(T)\leq n-b(T)+\rho(T)$.
\end{theorem}

This bound is sharp for generalized spiders (see Proposition \ref{Prop=}) and
for some caterpillars (Corollary \ref{Cor_Cat}). Combining Theorems
\ref{Thm_LB} and \ref{ThmMain} gives the following result.

\begin{theorem}
\label{Thm_Trees_both}For any tree $T$ of order $n$ that is not a path,%
\[
n-b(T)-|W_{\operatorname{int}}(T)|+\alpha(\operatorname{Int}(T))\leq
\alpha_{\operatorname{bn}}(T)\leq n-b(T)+\rho(T).
\]

\end{theorem}

When the tree $T$ has no internal vertices of degree $2$ (i.e., the only
vertices of degree $2$ lie on endpaths), we have the following corollary.

\begin{corollary}
\label{Wint=empty}If $T$ is a tree of order $n$ such that $b(T)\geq1$ and
$W_{\operatorname{int}}(T)=\varnothing$, then%
\[
n-b(T)+\alpha(T[R(T)])\leq\alpha_{\operatorname{bn}}(T)\leq n-b(T)+\rho(T).
\]

\end{corollary}

\noindent\textbf{Proof}.\hspace{0.1in}If $W_{\operatorname{int}}%
(T)=\varnothing$, then $\operatorname{Int}(T)$ is the subgraph of $T$ induced
by $R(T)$ and the result follows from Theorem \ref{Thm_Trees_both}%
.~$\blacksquare$

\bigskip

Note that the difference between the upper and lower bounds in Corollary
\ref{Wint=empty} equals $|R(T)|-\alpha(T[R(T)])$. The following question was
posed in \cite{MN3}:

\begin{question}
\label{Q1}\emph{\cite{MN3}}\hspace{0.1in}Can the upper bound in Theorem
\ref{ThmMain} be improved to $\alpha_{\operatorname{bn}}(T)\leq
|V(T)|-b(T)+\alpha(T[R(T)]$?
\end{question}

\begin{figure}[ptb]%
\centering
\includegraphics[
height=1.4814in,
width=3.9608in
]%
{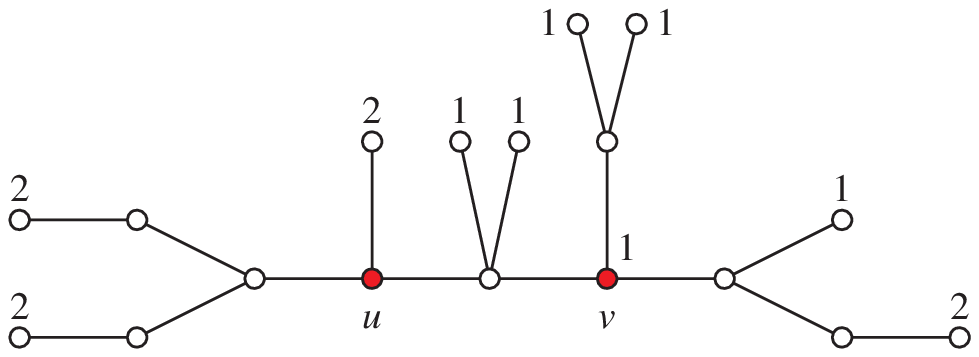}%
\caption{A tree $T$ of order $18$ such that $W_{\operatorname{int}%
}(T)=\varnothing$, $B_{0}(T)=\{v\},\ B_{1}(T)=\{u\}$, $R(T)=\{u,v\}$, and
$\alpha_{\operatorname{bn}}(T)=n-b(T)+\alpha(T[R(T)])=n-b(T)+\rho(T)=14.$}%
\label{Fig_equal}%
\end{figure}

If the answer to Question \ref{Q1} turns out to be yes, Corollary
\ref{Wint=empty} would give us an exact formula for $\alpha_{\operatorname{bn}%
}$ for trees with branch vertices but no internal degree $2$ vertices, namely
$\alpha_{\operatorname{bn}}(T)=|V(T)|-b(T)+\alpha(T[R(T)]$. As it stands,
Corollary \ref{Wint=empty} gives us an exact formula for a subclass of these
trees, an example of which is depicted in Figure \ref{Fig_equal}.%

\begin{corollary}
\label{Cor_RT_Ind}If $T$ is a tree of order $n$ such that $b(T)\geq1$,
$W_{\operatorname{int}}(T)=\varnothing$, and $R(T)$ is either empty or an
independent set, then%
\[
\alpha_{\operatorname{bn}}(T)=n-b(T)+\rho(T).
\]

\end{corollary}

\noindent\textbf{Proof}.\hspace{0.1in}If $R(T)$ is empty or independent, then
$\alpha(T[R(T)])=|R(T)|=\rho(T)$ and the result follows from Corollary
\ref{Wint=empty}.~$\blacksquare$

\bigskip

We next state a formula for $\alpha_{\operatorname{bn}}(T)$, where $T$ is a
tree with exactly two branch vertices; its proof can be found in Section
\ref{Thm_b(T)2_pf}. We use this result to show that there exist trees for
which $\alpha_{\operatorname{bn}}$ lies strictly between the bounds given in
Theorem \ref{Thm_Trees_both}. All of these trees support a positive answer to
Question \ref{Q1}.

\begin{theorem}
\label{b(T)2}{If $T$ is a tree of order $n$ such that $B(T)=\{b_{1},b_{2}\}$,
then }%
\[
\alpha_{\operatorname{bn}}{(T)=n-1-\min\{\lceil\tfrac{1}{2}d(b_{1}%
,b_{2})\rceil,\operatorname{loss}(b_{1}),\operatorname{loss}(b_{2})\}.}%
\]

\end{theorem}

\section{Known results}

\label{Sec_known}In this section we present known results that are of interest
or will be used later on. It is often useful to know when a bn-independent
broadcast $f$ is \emph{maximal bn-independent}, that is, there does not exist
a bn-independent broadcast $g$ such that $f<g$.

\begin{proposition}
\label{prop-max-bn}

\begin{enumerate}
\item[$(i)$] \emph{\cite[Proposition 2.2]{MN}}\hspace{0.1in}A bn-independent
broadcast $f$ on a graph $G$ is maximal bn-indepen-dent if and only if it is
dominating and either $V_{f}^{+}=\{v\}$ or $B_{f}(v)-\operatorname{PB}%
_{f}(v)\neq\varnothing$ for each $v\in V_{f}^{+}$.

\item[$(ii)$] \emph{\cite[Proposition 3.1]{MM}\hspace{0.1in}}Let $f$ be a
bn-independent broadcast on a connected graph $G$ such that $|V_{f}^{+}|\geq
2$. Then $f$ is maximal bn-independent if and only if each component of
$G-U_{f}^{E}$ contains at least two broadcasting vertices.
\end{enumerate}
\end{proposition}

Suppose $f$ is a bn-independent broadcast on $G$ and an edge $uv$ of $G$ is
covered by vertices $x,y\in V_{f}^{+}$. By the definition of covered,
$\{u,v\}\nsubseteq B_{f}(x)$ and $\{u,v\}\subseteq N_{f}(x)\cap N_{f}(y)$.
This violates the bn-independence of $f$. Hence we have the following observation.

\begin{observation}
\label{edge-covered}\emph{\cite[Observation 2.3]{MN}}\hspace{0.1in}If $f$ is a
bn-independent broadcast on a graph $G$, then each edge of $G$ is covered by
at most one vertex in $V_{f}^{+}$.
\end{observation}

The following bound on $\alpha_{\operatorname{bn}}(G)$ was proved in \cite{MN,
LindaD}.

\begin{proposition}
\label{Prop=}\emph{\cite[Corollaries 2.6, 2.7]{MN}}\hspace{0.1in}For any
connected graph $G$ of order $n$ and any spanning tree $T$ of $G$,
$\alpha_{\operatorname{bn}}(G)\leq\alpha_{\operatorname{bn}}(T)\leq n-1$.
Moreover, $\alpha_{\operatorname{bn}}(G)=n-1$ if and only if $G$ is a path or
a generalized spider.
\end{proposition}

Consider caterpillars $T$ such that $W_{\operatorname{int}}(T)=\varnothing$,
i.e., there are no vertices of degree $2$ between the first and last branch
vertices on the spine, and such that $B_{1}(T)$ is either empty or an
independent set. For such a caterpillar $T$, let $b_{1},...,b_{k}$ be the
branch vertices of $T$, labelled from left to right on the spine. Note that
$b_{1},b_{k}\in B_{\geq2}(T)$, i.e., $b_{1},b_{k}\notin R(T)$. If $v$ is a
branch vertex of a caterpillar $T$, then $v$ is a stem, hence $B_{0}%
(T)=\varnothing$ and $R(T)=B_{1}(T)$. The next result, which was obtained in
\cite{MN3} by explicitly defining a bn-independent broadcast of the
appropriate weight, now follows from Corollary \ref{Cor_RT_Ind}.

\begin{corollary}
\label{Cor_Cat}\emph{\cite[Corollary 5.1]{MN3}}\hspace{0.1in}If $T$ is a
caterpillar whose branch vertices induce a path $P=(b_{1},...,b_{k})$, and
$R(T)$ (i.e., the branch vertices among $b_{2},...,b_{k-1}$ that are adjacent
to exactly one leaf) is either empty or an independent set, then
$\alpha_{\operatorname{bn}}(T)=|V(T)|-b(T)+\rho(T)$.
\end{corollary}

To establish the upper bound in Theorem \ref{ThmMain}, several lemmas were
proved in \cite{MN3}. For simplicity we state those that we need here as
separate items of a single lemma.

\begin{lemma}
\label{Lem_hear1}$(i)\hspace{0.1in}$\emph{\cite[Lemma 3.1]{MN3}\hspace{0.1in}%
}For any tree $T$ and any $\alpha_{\operatorname{bn}}(T)$-broadcast $f$, no
leaf of $T$ hears $f$ from any non-leaf vertex.

\begin{enumerate}
\item[$(ii)$] \label{Lem_nonleaf1}\emph{\cite[Lemma 3.2]{MN3}\hspace{0.1in}%
}For any tree $T$ there exists an $\alpha_{\operatorname{bn}}(T)$-broadcast
$f$ such that $f(v)=1$ whenever $v\in V_{f}^{+}-L(T)$.

\item[$(iii)$] \label{Lem_leaf_branch}\emph{\cite[Lemma 3.4]{MN3}%
\hspace{0.1in}}Let $f$ be an $\alpha_{\operatorname{bn}}$-broadcast on a tree
$T$ such that a leaf $l$ $f$-dominates a branch vertex $w$. If $l^{\prime}$ is
a leaf in $L(w)$ that does not hear $f$ from $l$, then $l^{\prime}\in
V_{f}^{+}$, the $l^{\prime}-w$ path contains a vertex $b\in B_{f}(l)$, and
$f(l^{\prime})=d(b,l^{\prime})$.
\end{enumerate}
\end{lemma}

The final lemma demonstrates the importance of the branch vertices of a tree
$T$ in determining the broadcast values which may be assigned to its leaves.
Informally, Lemma \ref{Lem_branch} states that in an $\alpha
_{\operatorname{bn}}(T)$-broadcast, a leaf $l$ never overdominates a branch
vertex $b$ by exactly $2$. Either $l$ overdominates a branch vertex $b$ by
exactly $1$ and $b$ has no leaves except possibly $l$, or $l$ overdominates
$b$ by at least $3$ and has exactly one vertex in its boundary, this vertex
being not on a $b^{\prime}-l^{\prime}$ path for any $b^{\prime}\in B(T)$ and
$l^{\prime}\in L(b^{\prime})$. In addition to $b$, $l$ may overdominate any
number of branch vertices by $3$ or more as long as it also overdominates all
of their leaves. Although a more general result was proved in \cite[Lemma
3.6]{MN3}, we state Lemma \ref{Lem_branch} for broadcasts in which only leaves
broadcast with strength greater than $1$ as this is all that we need here.

\begin{lemma}
\label{Lem_branch}\emph{\cite[Corollary 3.7]{MN3}\hspace{0.1in}}Let $T$ be a
tree with $b(T)\geq2$ and $\mathcal{F}^{\prime}$ the set of all $\alpha
_{\operatorname{bn}}(T)$-broadcasts in which only leaves broadcast with
strength greater than $1$. Let $\mathcal{F}$ be the set of broadcasts in
$\mathcal{F}^{\prime}$ with the minimum number of overdominated branch
vertices. Then there exists a broadcast $f\in\mathcal{F}$ that satisfies the
following statement:

For any leaf $l$, let $X$ be the set of all branch vertices overdominated by
$l$. If $X\neq\varnothing$ and $v\in B_{f}(l)$, then $v$ is neither the leaf
nor an internal vertex on any endpath of $T$. Moreover,

\begin{enumerate}
\item[$(i)$] there exists $w\in X$ such that $f(l)=d(l,w)+1$, and either
$L(w)=\{l\}$ and $X=\{w\}$, or $L(w)=\varnothing$ and $f(l)\geq d(l,w^{\prime
})+3$ for all $w^{\prime}\in X-\{w\}$, or

\item[$(ii)$] $f(l)\geq d(l,w)+3$ for all $w\in X$ and $|B_{f}(l)|=1$.
\end{enumerate}
\end{lemma}

\section{Proof of Theorem \ref{Thm_LB}}

\label{Sec_lower_pf}We are now ready to prove Theorem \ref{Thm_LB}, which we
restate for convenience.

\bigskip

\noindent\textbf{Theorem \ref{Thm_LB}\hspace{0.1in}}\emph{If }$T$\emph{ is a
tree such that }$b(T)\geq1$\emph{, then}%
\[
\alpha_{\operatorname{bn}}(T)\geq n-b(T)-|W_{\operatorname{int}}%
(T)|+\alpha(\operatorname{Int}(T)).
\]

\begin{figure}[ptb]%
\centering
\includegraphics[
height=1.9372in,
width=6.3218in
]%
{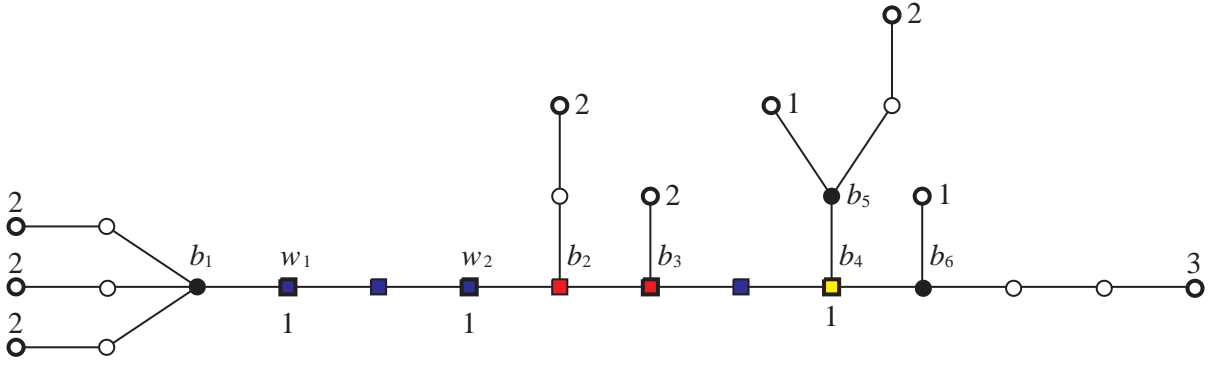}%
\caption{A tree $T$ with a bn-independent broadcast $f$ as described in the
proof of Theorem \ref{Thm_LB}. Note that $\sigma
(f)=n-b(T)-|W_{\operatorname{int}}(T)|+\alpha(\operatorname{Int}%
(T))=26-6-4+4=20<23=n-b(T)+\rho(T)=26-6+3$.}%
\label{Fig_lbT}%
\end{figure}

\noindent\textbf{Proof}.\hspace{0.1in}It is sufficient to define a
bn-independent broadcast $f$ on $T$ such that $\sigma
(f)=n-b(T)-|W_{\operatorname{int}}(T)|+\alpha(\operatorname{Int}(T))$. To do
this, let $X$ be a maximum independent set of $\operatorname{Int}(T)$ and let
$Y=B_{1}(T)-X$. The definition of $f$ below is illustrated in Figure
\ref{Fig_lbT}, where $X=\{w_{1},w_{2},b_{3},b_{4}\}$.

\begin{enumerate}
\item[$(i)$] For each vertex $b\in B_{\geq2}(T)\cup Y$ and each leaf $l\in
L(b)$, let $f(l)=d(b,l)$. 

Then $\bigcup_{l\in L(b)}N_{f}(l)=V(T_{b})$ and $\sum_{l\in L(b)}%
f(l)=|V(T_{b})|-1$.

\item[$(ii)$] For each vertex $b\in X\cap B_{1}(T)$ and the leaf $l\in L(b)$,
let $f(l)=d(b,l)+1$. 

Then $N_{f}(l)=V(T_{b})\cup N(b)$ and $f(l)=|V(T_{b})|$.

\item[$(iii)$] For each vertex $v\in X\cap(B_{0}(T)\cup W_{\operatorname{int}%
}(T))$, let $f(v)=1$.

\item[$(iv)$] Otherwise let $f(v)=0$.
\end{enumerate}

\noindent The weigth of $f$ is%
\begin{align*}
\sigma(f) &  =\sum_{b\in B_{\geq2}(T)\cup Y}\left(  |V(T_{b})|-1\right)
+\sum_{b\in X\cap B_{1}(T)}|V(T_{b})|+|X\cap(B_{0}(T)\cup
W_{\operatorname{int}}(T))|\\
&  =\sum_{b\in B_{\geq1}(T)}\left(  |V(T_{b})|-1\right)  +|X\cap
B_{1}(T)|+|X\cap(B_{0}(T)\cup W_{\operatorname{int}}(T))|\\
&  =\sum_{b\in B_{\geq1}(T)}\left(  |V(T_{b})|-1\right)  +|X|.
\end{align*}%
Since the expression $\sum_{b\in B_{\geq1}(T)}\left(  |V(T_{b})|-1\right)  $
counts all vertices in the subtrees $T_{b}$ except the branch vertex itself,
and $T_{b}=K_{1}$ if $b\in B_{0}(T)$,
\[
\sum_{b\in B_{\geq1}(T)}\left(  |V(T_{b})|-1\right)
=n-b(T)-|W_{\operatorname{int}}(T)|,
\]
from which we obtain%
\begin{align*}
\sigma(f) &  =n-b(T)-|W_{\operatorname{int}}(T)|+|X|\\
&  =n-b(T)-|W_{\operatorname{int}}(T)|+\alpha(\operatorname{Int}(T)).
\end{align*}

We show that $f$ is bn-independent. Suppose this is not the case. Then there
are two vertices $u,v\in V_{f}^{+}$ and an edge $xy\in E(T)$ such that the
broadcasts from $u$ and $v$ overlap on $xy$, that is, $\{x,y\}\subseteq
N_{f}(u)\cap N_{f}(v)$ but $\{x,y\}\nsubseteq B_{f}(u)\cap B_{f}(v)$. By the
definition of $f$, $u$ satisfies $(i)$, $(ii)$ or $(iii)$. We consider the
three cases separately and indicate the end of the proof of each case by a
solid diamond ($\blacklozenge$).\smallskip

\noindent\textbf{Case 1:}\hspace{0.1in}$u$ satisfies the conditions of $(i)$.
Then there exists a vertex $b\in B_{\geq1}(T)-X$ such that $u\in L(b)$ and
$f(u)=d(u,b)$. Therefore both $x$ and $y$ lie on the $u-b$ path in $T$ and at
least one of $x$ and $y$ is an internal vertex of this path. Since $v$
$f$-dominates $x$ and $y$, $f(v)\geq d(v,b)+1$. We show that this is not
possible. First suppose that $v\in L(b^{\prime})$ for some $b^{\prime}\in
B_{\geq1}(T)-X$. If $b^{\prime}=b$, then $f(v)=d(v,b)<d(v,b)+1$, and if
$b^{\prime}\neq b$, then $f(v)=d(v,b^{\prime})<d(v,b)$. Suppose next that
$v\in L(b^{\prime})$ for some $b^{\prime}\in X\cap B_{1}(T)$. Then $b\neq
b^{\prime}$. By the definition of $f$, $f(v)=d(v,b^{\prime})+1\leq
d(v,b)<d(v,b)+1$. Finally, suppose $v\in X\cap(B_{0}(T)\cup
W_{\operatorname{int}}(T))$. Then $v\neq b$ and $f(v)=1$, hence $f(v)\leq
d(v,b)<d(v,b)+1$. Each case produces a contradiction.~$\blacklozenge
$\smallskip

\noindent\textbf{Case 2:}\hspace{0.1in}$u$ satisfies the conditions of $(ii)$.
Then $u\in L(b)$, where $b\in X\cap B_{1}(T)$ and $f(u)=d(u,b)+1$, and at
least one of $x$ and $y$ lies on the $u-b$ path in $T$. Hence $f(v)\geq
d(v,b)$. We showed in Case 1 that $v$ does not satisfy $(i)$. Hence either
$v\in L(b^{\prime})$, where $b^{\prime}\in X\cap B_{1}(T)$ and
$f(v)=d(v,b^{\prime})+1$, or $f(v)=1$ and $v\in X\cap(B_{0}(T)\cup
W_{\operatorname{int}}(T))$. Since $f(v)\geq d(v,b)$, the former condition
implies that $b$ and $b^{\prime}$ are adjacent (note that $b\neq b^{\prime}$
because $u\neq v$ and $b,b^{\prime}\in B_{1}(T)$), while the latter condition
implies that $v$ and $b$ are adjacent. This contradicts the independence of
$X$.~$\blacklozenge$\smallskip

\noindent\textbf{Case 3:}\hspace{0.1in}$u,v\in X\cap(B_{0}(T)\cup
W_{\operatorname{int}}(T))$. In this case $f(u)=f(v)=1$, hence $d(u,v)=1$. But
$u,v\in X$, an independent set, and we again have a
contradiction.~$\blacklozenge$\smallskip

This exhausts all possible cases, hence $f$ is bn-independent. We deduce that
$\alpha_{\operatorname{bn}}(T)\geq\sigma(f)=n-b(T)-|W_{\operatorname{int}%
}(T)|+\alpha(\operatorname{Int}(T))$.~$\blacksquare$

\section{Proof of Theorem \ref{b(T)2}}

\label{Thm_b(T)2_pf}In this section we fist prove Theorem \ref{b(T)2}, then we
use it to show that there exist trees for which $\alpha_{\operatorname{bn}}$
lies strictly between the bounds in Theorem \ref{Thm_Trees_both}. We again
restate the theorem for convenience.

\medskip

\noindent\textbf{Theorem \ref{b(T)2}\hspace{0.1in}}\emph{If }$T$\emph{ is a
tree of order }$n$\emph{ such that }$B(T)=\{b_{1},b_{2}\}$\emph{, then}{ }%
\[
\alpha_{\operatorname{bn}}{(T)=n-1-\min\{\lceil\tfrac{1}{2}d(b_{1}%
,b_{2})\rceil,\operatorname{loss}(b_{1}),\operatorname{loss}(b_{2})\}.}%
\]

\bigskip

\noindent\textbf{Proof.\hspace{0.1in}}Let $\mathcal{F}^{\prime}$ be the set of
all $\alpha_{\operatorname{bn}}(T)$-broadcasts in which only leaves broadcast
with strength greater than $1$. By Lemma \ref{Lem_nonleaf1}$(ii)$,
$\mathcal{F}^{\prime}\neq\varnothing$. Let $\mathcal{F}$ be the set of
broadcasts in $\mathcal{F}^{\prime}$ with the minimum number of overdominated
branch vertices and choose a broadcast $f\in\mathcal{F}$ such that Lemma
\ref{Lem_branch} applies. By Proposition \ref{prop-max-bn}$(i)$, all vertices
of $T$ are $f$-dominated. By the choice of $\mathcal{F}^{\prime}$, each
non-leaf vertex of $T$ is $f$-dominated by a leaf or by a broadcast of
strength $1$, and (by Lemma \ref{Lem_nonleaf1}$(i)$) each leaf is
$f$-dominated by a leaf.

Let $P=(b_{1}=v_{0},v_{1},...v_{k}=b_{2})$ be the $b_{1}-b_{2}$ path in $T$.
Notice that $d(b_{1},b_{2})=k$. Suppose $k=1$. Then there are no internal
vertices on $P$ and $R(T)=\varnothing$. Also,
\[
{\min\{\lceil\tfrac{1}{2}d(b_{1},b_{2})\rceil,\operatorname{loss}%
(b_{1}),\operatorname{loss}(b_{2})\}=\lceil\tfrac{1}{2}d(b_{1},b_{2})\rceil
=1}.
\]
By Corollary \ref{Cor_RT_Ind},
\[
\alpha_{\operatorname{bn}}(T)=n-b(T)=n-2=n-1-{\min\{\lceil\tfrac{1}{2}%
d(b_{1},b_{2})\rceil,\operatorname{loss}(b_{1}),\operatorname{loss}(b_{2})\}}%
\]
as required. We therefore assume that $k\geq2$ and examine the way in which
the internal vertices of $P$ are dominated. Let $P_{\operatorname{int}}%
:v_{1},...,v_{k-1}$ be the path induced by the internal vertices of $P$. There
are three possibilities for dominating $P_{\operatorname{int}}$. We indicate
the end of the proof of each case by a solid diamond ($\blacklozenge$).
\smallskip

\noindent\textbf{Case 1:\hspace{0.1in}}No vertex on $P_{\operatorname{int}}$
is $f$-dominated by a leaf. We show that $\sigma(f)=n-1-{\lceil\tfrac{1}%
{2}d(b_{1},b_{2})\rceil}$.

The vertices on $P_{\operatorname{int}}$ are $f$-dominated by a set $D$ of
vertices on $P_{\operatorname{int}}$. By Lemma \ref{Lem_hear1}$(iii)$,
$f(x)=d(x,b_{1})$ for every leaf $x\in L(b_{1})$. Similarly, for every leaf
$y\in L(b_{2})$, we have $f(y)=d(y,b_{2})$. For $i\in\{1,2\}$, let $f_{i}$ be
the restriction of $f$ to $T_{b_{i}}$ (the subtree induced by the vertices on
the $b_{i}-l$ paths for all $l\in L(b_{i})$). Notice that $\sigma
(f_{i})=|V(T_{i})|-1$ for $i=1,2$; by Proposition \ref{Prop=}, $f_{i}$ is an
$\alpha_{\operatorname{bn}}(T_{b_{i}})$-broadcast. By the choice of
$f\in\mathcal{F}^{\prime}$ and Lemma \ref{Lem_nonleaf1}$(ii)$, $f(x)=1$ for
each $x\in D$. Since $f$ is a maximum bn-independent broadcast, $D$ is a
maximum independent set of $P_{\operatorname{int}}$, hence $|D|=\alpha
(P_{\operatorname{int}})=\left\lceil \frac{1}{2}(k-1)\right\rceil $. It
follows that
\begin{align*}
\sigma(f)  &  =|V(T_{b_{1}})|-1+|V(T_{b_{2}})|-1+\left\lceil \frac{1}%
{2}(k-1)\right\rceil \\
&  =n-2-\left\lfloor \frac{1}{2}(k-1)\right\rfloor \\
&  =n-1-{\lceil\tfrac{1}{2}d(b_{1},b_{2})\rceil,}%
\end{align*}
which is what we wanted to show.~$\blacklozenge$

\noindent\textbf{Case 2:\hspace{0.1in}}A leaf $x\in L(b_{1})$ $f$-dominates
some or all of the vertices on $P_{\operatorname{int}}$ and thus overdominates
$b_{1}$. We show that $\sigma(f)=n-1-{\operatorname{loss}(b_{1})}$.

Since $|L(b_{1})|\geq2$, Lemma \ref{Lem_branch}$(i)$ does not apply to $x$.
Thus Lemma \ref{Lem_branch}$(ii)$ applies and $B_{f}(x)=\{v\}$, where $v$ is a
vertex on $P-b_{1}$. Hence $x$ dominates $T_{b_{1}}$ and overdominates
$L(b_{1})$ (otherwise a vertex on a $b_{1}-l$ endpath belongs to
$B_{f}(x)-\{v\}$, which is impossible).\label{here}

Let $v_{j}$ be the vertex on the path $P$ such that $B_{f}(x)=\{v_{j}\}$. We
state and prove two claims, and indicate the end of the proof of each claim by
an open diamond ($\lozenge$).

\begin{claim}
\label{Cl1}$B_{f}(x)=\{b_{2}\}$, i.e., $v_{j}=v_{k}=b_{2}$, hence
$f(x)=d(x,b_{2})$.
\end{claim}

\noindent\textbf{Proof of Claim \ref{Cl1}\hspace{0.1in}}Suppose first that
$j=k-1$. Then $b_{2}$ does not hear $f$ from $x$; hence some vertex $y\in
V(T_{b_{2}})$ broadcasts to $b_{2}$. Suppose $y\in L(b_{2})$. Since
$|L(b_{2})|\geq2$, Lemma \ref{Lem_branch} implies that $y$ does not
overdominate $b_{2}$ by exactly $1$. Since $f$ is bn-independent, $y$ does not
broadcast to $v_{k-2}$. Therefore $f(y)=d(y,b_{2})$ and $b_{2}\in B_{f}(y)$.
But then $\{v_{k-1}\}=\operatorname{PB}_{f}(x)=B_{f}(x)$, and by Proposition
\ref{prop-max-bn}$(i)$, $f$ is not maximal bn-independent. Suppose $y=b_{2}$.
By Lemma \ref{Lem_nonleaf1}$(ii)$, $f(b_{2})=1$. For each $l\in L(b_{2})$, let
$Q_{l}^{\prime}$ be the $b_{2}-l$ endpath and $Q_{l}=Q_{l}^{\prime}-b_{2}$.
Since $f(b_{2})=1$, $\sum_{v\in V(Q_{l})}f(v)\leq d(l,b_{2})-1$. Define the
broadcast $g_{1}$ by
\[
g_{1}(u)=\left\{
\begin{tabular}
[c]{rl}%
$0$ & if $u=b_{2}$,\\
$0$ & if $u$ is a degree $2$ vertex on $Q_{l}$ for $l\in L(b_{2})$,\\
$d(u,b_{2})$ & if $u\in L(b_{2})$,\\
$f(u)$ & otherwise.
\end{tabular}
\right.
\]
Then $g_{1}$ is bn-independent and, since $|L(b_{2})|\geq2$, $\sigma
(g_{1})>\sigma(f)$, which is impossible. We obtain a similar contradiction if
$y$ is an internal vertex on $Q_{l}$ for some $l\in L(B_{2})$. Therefore
$j\neq k-1$.

Suppose next that $j<k-1$. By our choice of $f$ and because $f$ is dominating,
$v_{j+1}$ is $f$-dominated either by $t\in\{v_{j+1},v_{j+2}\}$, where
$f(t)=1$, or by a leaf $l\in L(b_{2})$. In the former case, the broadcast
$g_{2}$ with $g_{2}(x)=d(x_{2},t)+1$, $g_{2}(t)=0$ and $g_{2}(u)=f(u)$
otherwise, has greater weight than $f$, and $N_{g_{2}}(x)\subseteq
N_{f}(t)\cup N_{f}(x)$. Hence $g_{2}$ is bn-independent and violates the
maximality of $f$. In the latter case, Lemma \ref{Lem_branch}$(ii)$ applies to
$l$, hence $l$ overdominates $L(b_{2})$ and $f(l)=d(l,v_{j})$, $f(x)=d(x,v_{j}%
)$ and $f(u)=0$ otherwise.

Define a new broadcast $g_{3}$ by $g_{3}(u)=d(u,b_{2})$ if $u=x$ or $u\in
L(b_{2})$, and $g(u)=0$ otherwise. For all $\{u,u^{\prime}\}\in V_{g_{3}}^{+}$
such that $u\neq u^{\prime}$, $N_{g_{3}}(u)\cap N_{g_{3}}(u^{\prime}%
)=\{b_{2}\}=B_{g_{3}}(u)\cap B_{g_{3}}(u^{\prime})$, hence $g_{3}$ is
bn-independent. Since $g_{3}(x)+g_{3}(l)=f(x)+f(y)$ and $|L(b_{2})|\geq2$,
$\sigma(g_{3})>\sigma(f)$. Thus, $g_{3}$ violates the maximality of $f$. We
deduce that $f(x)=d(x,b_{2})$, that is, $B_{f}(x)=\{b_{2}\}$.~$\lozenge
$\smallskip

Proceeding with the proof of Case 2, we point out that $b_{2}$ is the only
vertex of $T_{b_{2}}$ that is $f$-dominated by $x$. As in Case 1, let $f_{2}$
be the restriction of $f$ to $T_{b_{2}}$. Then $f_{2}$ dominates all vertices
of $T_{b_{2}}$ except possibly $b_{2}$. By Proposition \ref{Prop=},
$\alpha_{\operatorname{bn}}(T_{b_{2}})=|V(T_{b_{2}})|-1$ and thus
$\sigma(f_{2})\leq|V(T_{b_{2}})|-1$. Moreover, the broadcast $h$ on $T_{b_{2}%
}$ defined by $h(l)=d(l,b_{2})$ for each $l\in L(b_{2})$ is a bn-independent
broadcast of weight $\sigma(h)=|V(T_{b_{2}})|-1$ such that $(f-f_{2})\cup h$
is a bn-independent broadcast on $T$ with weight at least the weight of $f$.
The maximality of $f$ now implies that
\begin{equation}
\sigma(f_{2})=|V(T_{b_{2}})|-1. \label{eq_f2a}%
\end{equation}

We state and prove our second claim.

\begin{claim}
\label{Cl2}$d(x,b_{1})=\max(b_{1})$.
\end{claim}

\noindent\textbf{Proof of Claim \ref{Cl2}\hspace{0.1in}}Suppose to the
contrary that $d(x,b_{1})<\max(b_{1})$ and let $x^{\prime}$ be a leaf such
that $d(x^{\prime},b_{1})=\max(b_{1})$. Create a new broadcast $g_{4}$ with
$g_{4}(x)=0$, $g_{4}(x^{\prime})=d(x^{\prime},b_{2})$ and $g_{4}(u)=g_{4}(u)$
otherwise. Since $d(x,b_{1})<d(x^{\prime},b_{1})$, $\sigma(g_{4})>\sigma
(g_{4})$. Notice that $N_{g_{4}}(x^{\prime})=N_{g_{4}}(x)$ and all other
boundaries are unchanged. Hence $g_{4}$ is bn-independent and contradicts the
maximality of $f$.~$\lozenge$\smallskip

Recall that by definition, {$\operatorname{loss}(b_{1})=\operatorname{sum}%
(b_{1})-\max(b_{1})=$}$\sum_{l\in L(b_{1})}${$d_{T}(b_{1},l)$}$-\max(b_{1})${,
hence }$\max(b_{1})=|V(T_{b_{1}})|-1-\operatorname{loss}(b_{1})$. By
(\ref{eq_f2a}) and Claims \ref{Cl1} and \ref{Cl2} we now have that%
\begin{align*}
\sigma(f)  &  =f(x)+\sigma(f_{2})=d(x,b_{1})+d(b_{1},b_{2})+\sigma(f_{2})\\
&  =\max(b_{1})+d(b_{1},b_{2})+\sigma(f_{2})\\
&  =|V(T_{b_{1}})|-1-\operatorname{loss}(b_{1})+k+|V(T_{b_{2}})|-1.
\end{align*}
Since $V(T)$ is the disjoint union of $V(T_{b_{1}})-\{b_{1}\}$, $V(P)$ and
$V(T_{b_{2}})-\{b_{2}\}$, and $P$ is a path of order $k+1$, it follows that
$\sigma(f)=n-1-{\operatorname{loss}(b_{1})}$.~$\blacklozenge$\smallskip

\noindent\textbf{Case 3:\hspace{0.1in}}A leaf $x\in L(b_{2})$ $f$-dominates
some or all of the vertices on $P_{\operatorname{int}}$ and thus overdominates
$b_{2}$. Exactly as in Case 2, $\sigma(f)=n-1-{\operatorname{loss}(b_{2})}%
$.~$\blacklozenge$\smallskip

Since $f$ is a bn-independent broadcast of maximum weight, it follows from
Cases 1, 2 and 3 that $\alpha_{\operatorname{bn}}(T)=\sigma(f)={n-1-\min
\{\lceil\tfrac{1}{2}d(b_{1},b_{2})\rceil,\operatorname{loss}(b_{1}%
),\operatorname{loss}(b_{2})\}}$.~$\blacksquare$

\bigskip

\begin{figure}[ptb]%
\centering
\includegraphics[
height=0.6944in,
width=4.6216in
]%
{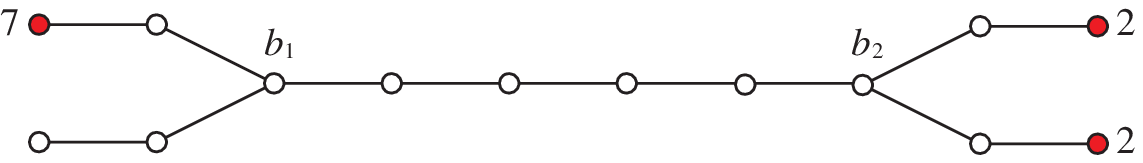}%
\caption{A tree $T$ with $b(T)=2$, $R(T)=\varnothing$, $|W_{\operatorname{int}%
}(T)|=4$ and $\alpha(\operatorname{Int}(T))=2$. Hence
$n-b(T)-|W_{\operatorname{int}}(T)|+\alpha(\operatorname{Int}(T))=10<\alpha
_{\operatorname{bn}}(T)=11<n-2=12.$}%
\label{Fig_lb4}%
\end{figure}

Theorem \ref{b(T)2} allows us to determine a special class of trees $T$ for
which $\alpha_{\operatorname{bn}}(T)$ lies strictly in between our upper and
lower bounds, an example of which is given in Figure \ref{Fig_lb4}.%

\begin{corollary}
\label{between}{Let $T$ be a tree with branch set $B(T)=\{b_{1},b_{2}\}$ such
that $2\leq\operatorname{loss}(b_{1})\leq\operatorname{loss}(b_{2}%
)<\lceil\frac{d(b_{1},b_{2})}{2}\rceil$. Then
\[
n-b(T)-|W_{\operatorname{int}}(T)|+\alpha(\operatorname{Int}(T))<\alpha
_{\operatorname{bn}}(T)<n-b(T)-\rho(T).
\]
}
\end{corollary}

\noindent\textbf{Proof.\hspace{0.1in}}For a tree $T$ with branch number
$b(T)=2$, $B_{0}(T)\cup B_{1}(T)=\varnothing$, hence $\rho(T)=0$. Therefore
$n-b(T)-\rho(T)=n-2$. Also, $\operatorname{Int}(T)$ is the subgraph of $T$
induced by $W_{\operatorname{int}}(T)$, which, in this case, is the path
$P_{\operatorname{int}}$ of order $k-1=d(b_{1},b_{2})-1$ as described above in
the proof of Theorem \ref{b(T)2}. Hence
\begin{align*}
n-b(T)-|W_{\operatorname{int}}(T)|+\alpha(\operatorname{Int}(T))  &
=n-2-(k-1)+\left\lceil \frac{k-1}{2}\right\rceil \\
& =n-2-\left\lfloor \frac{k-1}{2}\right\rfloor =n-2-\left\lceil {\frac
{d(b_{1},b_{2})}{2}}\right\rceil {.}%
\end{align*}
By Theorem \ref{b(T)2} and the choice of $T$,%
\[
n-2-\left\lceil {\frac{d(b_{1},b_{2})}{2}}\right\rceil
<n-1-\operatorname{loss}(b_{1})=\alpha_{\operatorname{bn}}%
(T)<n-2.~\blacksquare
\]
We close this section by remarking that strict inequality in the upper bound
suggested in Question \ref{Q1} holds for the trees described in the statement
of Corollary \ref{between} as well. 

\section{Open problems}

\label{Sec_Open}In this final section we mention some open problems for future
research, beginning by repeating Question \ref{Q1} for the sake of completeness.

\bigskip

\noindent\textbf{Question 1\hspace{0.1in}}\cite{MN3}\hspace{0.1in}\emph{Can
the upper bound in Theorem \ref{ThmMain} be improved to }$\alpha
_{\operatorname{bn}}(T)\leq|V(T)|-b(T)+\alpha(T[R(T)]$\emph{?}

\bigskip

It is well known that $\alpha(G)\leq n-\delta(G)$ and, when $G$ is connected,
$\operatorname{diam}(G)\leq n-\delta(G)$.

\begin{question}
\label{Q_bn_upperbound}\emph{\cite{MN2}}\hspace{0.1in}Is it true that
$\alpha_{\operatorname{bn}}(G)\leq n-\delta(G)$ for all graphs $G$?
\end{question}

\begin{problem}
Characterize trees $T$ such that 

\begin{enumerate}
\item[$(i)$] \emph{\cite{MN3}}\hspace{0.1in}$\alpha_{\operatorname{bn}%
}(T)=|V(T)|-b(T)+\rho(T)$, i.e., equality holds in the upper bound in
Theorem~\ref{Thm_Trees_both},

\item[$(ii)$] $\alpha_{\operatorname{bn}}(T)=n-b(T)-|W_{\operatorname{int}%
}(T)|+\alpha(\operatorname{Int}(T))$, i.e., equality holds in the lower bound
in Theorem~\ref{Thm_Trees_both},

\item[$(iii)$] $\alpha_{\operatorname{bn}}(T)=|V(T)|-b(T)+\rho
(T)=n-b(T)-|W_{\operatorname{int}}(T)|+\alpha(\operatorname{Int}(T))$, i.e.,
the bounds in Theorem~\ref{Thm_Trees_both} coincide.
\end{enumerate}
\end{problem}

By definition, $\alpha(G)\leq\alpha_{\operatorname{bn}}(G)\leq\alpha_{h}(G)$
for all graphs $G$. Mynhardt and Neilson \cite{MN} showed that $\alpha
(G)=\alpha_{\operatorname{bn}}(G)$ when $G$ is a $2$-connected bipartite
graph. The spiders $\operatorname{Sp}(2^{k})$, for which $\alpha
_{\operatorname{bn}}(\operatorname{Sp}(2^{k}))=2k$ and $\alpha
(\operatorname{Sp}(2^{k}))=k+1$, demonstrate that this result does not hold
for trees.

\begin{problem}
Investigate the ratio $\alpha_{\operatorname{bn}}(G)/\alpha(G)$ when $G$ is
$(i)$ a tree, $(ii)$ a cyclic bipartite graph of connectivity $1$, $(iii)$
general graphs.
\end{problem}

\begin{problem}
\emph{\cite{MN3}}\hspace{0.1in}Characterize trees $T$ such that $(i)$
$\alpha_{\operatorname{bn}}(T)=\alpha_{h}(T)$ or $(ii)$ $\alpha
_{\operatorname{bn}}(T)=\alpha(T)$.
\end{problem}

\begin{problem}
\emph{\cite{MN3}}\hspace{0.1in}Characterize caterpillars $T$ such that $(i)$
$\alpha_{\operatorname{bn}}(T)=\alpha_{h}(T)$ or $(ii)$ $\left.
\alpha_{\operatorname{bn}}(T)=\alpha(T).\right.  $
\end{problem}

\begin{problem}
\emph{\cite{MN3}}\hspace{0.1in}Determine $\alpha_{\operatorname{bn}}(T)$ for
all caterpillars $T$.
\end{problem}

\smallskip

\noindent\textbf{Acknowledgement\hspace{0.1in}}We acknowledge the support of
the Natural Sciences and Engineering Research Council of Canada (NSERC), PIN 253271.

\noindent Cette recherche a \'{e}t\'{e} financ\'{e}e par le Conseil de
recherches en sciences naturelles et en g\'{e}nie du Canada (CRSNG), PIN
253271.%
\begin{center}
\includegraphics[
natheight=0.773100in,
natwidth=1.599900in,
height=0.4151in,
width=0.8276in
]%
{../../../NSERC2014_2019/NSERC_DIGITAL_BW/NSERC_DIGITAL_BW/NSERC_BLACK.jpg}%
\end{center}

\end{document}